\documentclass[a4paper,10pt]{article}

\usepackage{amstext}
\usepackage{amsmath}
\usepackage{ntheorem}
\usepackage{ifthen}
\usepackage{amscd}
\usepackage{amssymb}
\usepackage[english,french]{babel}
\usepackage[T1]{fontenc}
\usepackage[latin1]{inputenc}
\usepackage[all]{xy}
\usepackage{graphicx}
\usepackage{enumerate}

\newenvironment{demo}[1][]{\ifthenelse{\equal{#1}{}}{\noindent \textbf{Démonstration :\\ \indent}}{\noindent \textbf{Démonstration #1 :\\ \indent}}}{$\square$\\}

\theoremstyle{break}
\newtheorem{theo}{Théorème}[section]
\newtheorem{conj}{Conjecture}[section]

\newtheorem{prop}{Proposition}[section]
\newtheorem{defi}{Définition}[section]
\theorembodyfont{\slshape}
\newtheorem{rem}{Remarque}[section]
\newtheorem{qt}{Question}[section]
\theorembodyfont{\upshape}
\newtheorem{ex}{Exemple}[section]

\newcommand{\RR}{\mathbb R}

\newcommand{\CC}{\mathbb C}
\newcommand{\QQ}{\mathbb Q}
\newcommand{\ZZ}{\mathbb Z}

\newcommand{\PP}{\mathbb P}
\newcommand{\GG}{\Gamma}
\newcommand{\adg}{\textsc{adg}}

\newcommand{\lie}{\mathcal{L}}

\newcommand{\gf}[1]{\pi_1(#1)}
\newcommand{\To}{\longrightarrow}

\newcommand{\wtx}{\widetilde{X}}

\title{Remarques sur la cohomologie des groupes kählériens nilpotents}
\author{Benoît \textsc{Claudon}}

\begin{document}

\maketitle

\begin{abstract}
Dans cette note, nous montrons que la cohomologie des groupes kählériens (virtuellement) nilpotents portent une structure de Hodge mixte naturelle, les morphismes de Hopf devenant des morphismes de structures de Hodge mixtes. Nous illustrons ce phénomène sur les exemples connus de groupes kählériens nilpotents (non abéliens).
\end{abstract}
\vspace*{0.5cm}
\begin{center}
{\LARGE Remarks on the cohomology of K\"{a}hler groups}
\end{center}
\selectlanguage{english}
\begin{abstract}
In this note, we show that the cohomology groups of (virtually) nilpotent K\"{a}hler groups are naturally endowed with a mixed Hodge structure. These structures make the Hopf morphisms into mixed Hodge structures morphisms. We illustrate this fact with the study of known examples of non-abelian nilpotent K\"{a}hler groups.
\end{abstract}
\selectlanguage{french}

\section{Introduction}

Soit $X$ une variété kählérienne compacte de groupe fondamental $\Gamma=\gf{X}$. La cohomologie de ces deux objets est relié par des morphismes naturels :
$$H^k(\Gamma,\CC)\To H^k(X,\CC),$$
dont l'existence est due à Hopf. Pour les petits degrés, on peut bien sûr être plus précis ; en degré 0 et 1, ces morphismes sont des isomorphismes. En degré 2, on dispose de plus d'une suite exacte courte :
\begin{displaymath}\label{hopf}
0\To H^2(\GG,\CC)\To H^2(X,\CC)\To H^2(\wtx,\CC),
\end{displaymath}
où $\wtx$ est le revêtement universel de $X$ et le morphisme $H^2(X,\CC)\To H^2(\wtx,\CC)$ étant celui induit par la projection naturelle $\wtx\To X$. En particulier, $H^2(\GG,\CC)$ est de dimension finie.

Il est communément admis que la cohomologie des groupes kählériens devrait se comporter d'une façon similaire à celle des variétés. En particulier, la conjecture suivante est attribuée à Carlson et Toledo (relayée notamment par Koll\'ar \cite{K95}) :
\begin{conj}\label{carlson-toledo}
Le groupe $H^2(\GG,\CC)$ est toujours non-nul (pour $\GG$ un groupe kählérien \emph{infini}).
\end{conj}
Cependant, on ne sait pas trop qu'elle devrait être la forme des énoncés en degré plus élevé ; en effet, contrairement à celle des variétés compactes, la cohomologie des groupes n'est pas nécessairement de dimension finie en degré $\geq 3$ (le premier exemple de groupe de présentation finie dont la cohomologie est de dimension infinie est dû à Stallings\footnote{le groupe défini par la présentation
$$\langle a,b,c,x,y\vert\,[x,a],\,[x,b],\,[y,a],\,[y,b],\,[a^{-1}x,c],\,[a^{-1}y,c],\,[b^{-1}a,c]\rangle$$
a un groupe de cohomologie de dimension inifinie en degré 3.} \cite{St63} ; pour les exemples kählériens, voir \cite{Dim09}).\\

Indépendamment de savoir si la réponse à la conjecture \ref{carlson-toledo} est affirmative, la question suivante est assez naturelle :
\begin{qt}\label{question}
le sous-espace vectoriel $H^2(\GG,\CC)$ est-il une sous-structure de Hodge de $H^2(X,\CC)$ ?
\end{qt}
En d'autres termes, on doit vérifier l'égalité :
\begin{align*}
H^2(\GG,\CC)=\left(H^2(\GG,\CC)\cap H^{2,0}(X,\CC)\right)&\oplus\left(H^2(\GG,\CC)\cap H^{1,1}(X,\CC)\right)\\
&\oplus\left(H^2(\GG,\CC)\cap H^{0,2}(X,\CC)\right).
\end{align*}

Dans sa plus grande généralité, la question \ref{question} semble hors de portée des techniques actuelles (ou du moins nécessiter une idée nouvelle). En revanche, dans le cas des groupes nilpotents (voir la section \ref{exemple-fred} pour les exemples de groupes kählériens nilpotents), nous allons constater que la réponse est affirmative.

\begin{theo}\label{th1}
Soit $X$ une variété kählérienne compacte de groupe fondamental $\GG$ virtuellement nilpotents. Les groupes de cohomologie $H^k(\GG,\CC)$ (de dimension finie) sont naturellement munis de structures de Hodge mixtes (fonctorielles). De plus, les morphismes naturels
$$H^k(\GG,\CC)\To H^k(X,\CC)$$
sont des morphismes de \textsc{\textsc{shm}}.
\end{theo}

\noindent En degré 2, on peut même être plus précis.
\begin{theo}\label{th2}
Sous les mêmes hypothèses que ci-dessus, la \textsc{\textsc{shm}} sur $H^2(\GG,\CC)$ est pure (de poids 2) et on a :
\begin{align*}
H^2(\GG,\CC)&=\mathrm{Im\left(H^1(X,\CC)\wedge H^1(X,\CC)\To H^2(X,\CC)\right)}\\
&=\mathrm{Im\left(H^2(\mathrm{Alb}(X),\CC)\stackrel{\alpha^*}{\To} H^2(X,\CC)\right)}.
\end{align*}
En particulier, la conjecture \ref{carlson-toledo} est vraie pour les groupes kählériens nilpotents.
\end{theo}

En utilisant \cite{D}, il suffit par exemple de supposer que le groupe fondamental de $X$ est résoluble.
\begin{theo}\label{delzant}
Si le groupe fondamental d'une variété kählérienne compacte est résoluble, alors il est virtuellement nilpotent.
\end{theo}

\section{Rappels}

\subsection{Groupes nilpotents}

Soit $G$ un groupe (de type fini) et $C^i(G)$ sa suite centrale descendante définie par $C^1(G)=G$ et $C^{i+1}(G)=[C^i(G),G]$ pour $i\geq1$. On notera $G_i=G/C^{i+1}(G)$ les quotients (nilpotents) correspondants. Les éléments d'ordre fini de $G_i$ forment un sous-groupe fini caractéristique noté $\mathrm{Tor}(G_i)$ et $G_i^*=G_i/\mathrm{Tor}(G_i)$ est donc un groupe nilpotent sans torsion. On peut donc lui appliquer appliquer la proposition suivante.
\begin{prop}[Mal\v{c}ev, \cite{M49}]\label{completion malcev}
Soit $N$ un groupe de type fini, nilpotent et sans torsion. Il existe un unique groupe de Lie nilpotent (défini sur $\QQ$) et simplement connexe $N_\RR$ et une injection $N\hookrightarrow N_\RR$ qui réalise $N$ comme un réseau cocompact de $N_\RR$. On notera $\lie(N)$ l'algèbre de Lie de $N_\RR$.
\end{prop}
Cette proposition montre qu'on peut associer au groupe $G$ une tour d'extensions centrales :
$$\dots\To \lie_{i+1}(G)\To \lie_i(G)\To\dots\To \lie_1(G)\To 0,$$
où l'on a noté $\lie_i(G):=\lie(G_i^*)$. La limite projective de cette suite est notée
$$\lie(G):=\varprojlim \lie_i(G)$$
et est appelée la complétion de Mal\v{c}ev de $G$ (si $G$ est nilpotent, la suite est finie et $\lie(G)=\lie(G/\mathrm{Tor}(G))$ est une algèbre de Lie nilpotente de dimension finie).

\subsection{Cohomologie des groupes}

Soit $G$ un groupe et $M$ un $G$-module (on sera surtout concerné par le cas du module trivial). L'assignation $F:M\mapsto M^G$ qui a un $G$ module associe le sous-module de ses éléments $G$-invariants est un foncteur de la catégorie des $G$-modules vers celle des groupes abéliens, qui est de plus exact à gauche. On définit alors la cohomologie de $G$ à valeurs dans $M$ comme le foncteur dérivé de $F$ :
$$\forall\, k\ge 0,\,H^k(G,M)=\mathrm{R}^kF(M).$$
C'est aussi la cohomologie du complexe $(\mathcal{C}^{\bullet}(G,M),d)$ où $\mathcal{C}^{k}(G,M)$ est constituée des applications de $G^k$ dans $M$ (par convention, $\mathcal{C}^{0}(G,M)=M$) et la différentielle étant donnée par :
\begin{align*}
df(g_1,\dots,g_{k+1})=g_1\cdot f(g_2,\dots,g_{k+1})+&\sum_{j=1}^k (-1)^jf(g_1,\dots,g_{j-1},g_jg_{j+1},\dots)\\
&+(-1)^{k+1}f(g_1,\dots,g_k).
\end{align*}

Dans le cas des groupes nilpotents sans torsion, la cohomologie se calcule facilement grâce à la complétion de Mal\v{c}ev.
\begin{theo}[K. Nomizu \cite{N54}, voir aussi \cite{R72}]\label{cohomologie nilp}
Si $G$ est un groupe nilpotent (de type fini) sans torsion d'algèbre de Lie $\lie(G)$, on dispose des isomorphismes suivants :
$$\forall k\ge0,\,H^k(G,\CC)\simeq H^k(\lie(G),\CC).$$
En particulier, $H^k(G,\CC)$ est de dimension finie pour tout $k\ge 0$.
\end{theo}
\begin{rem}\label{coho lie}
La cohomologie d'une algèbre de Lie $(\lie,[\cdot,\cdot])$ est celle du complexe des formes alternées $\Lambda^{\bullet}\lie^*$, la différentielle étant défini comme l'action duale du crochet de Lie $[\cdot,\cdot]:\lie\wedge\lie\To \lie$.
\end{rem}

Pour finir, mentionnons un outil très utile en cohomologie des groupes : l'opération de transfert. Soit $H\leq G$ un sous-groupe \emph{d'indice fini} de $G$. L'inclusion $H\stackrel{i}{\hookrightarrow} G$ induit un morphisme $i^*:H^{\bullet}(G,M)\To H^{\bullet}(H,M)$ mais, fait remarquable, il existe aussi un morphisme\footnote{la notation $V$ provient de l'allemand \emph{Verlagerung}.} allant dans la direction opposée :
$$V_{H\to G}=V:H^{\bullet}(H,M)\To H^{\bullet}(G,M)$$
et qui vérifie : $V\circ i^*=[G:H]\mathrm{Id}$. On a donc :
\begin{prop}[voir prop. 10.4, p. 85 \cite{B87}]\label{transfert}
Si la multiplication par $[G:H]$ est un automorphisme de $M$, l'application $i^*:H^{\bullet}(G,M)\To H^{\bullet}(H,M)$ est injective. Si de plus $H$ est un sous-groupe \emph{normal} de $G$, la cohomologie de $G$ s'identifie à la partie invariante sous l'action de $G/H$ de la cohomologie de $H$ :
$$H^{\bullet}(G,\CC)\stackrel{\sim}{\To} H^{\bullet}(H,\CC)^{G/H}.$$
\end{prop}
\begin{rem}\label{rm-transfert}
L'existence de l'application $V:G_{ab}\To H_{ab}$ avait d'abord été observée par Schur ; le transfert fut ensuite généralisé aux autres groupes de cohomologie par Eckmann \cite{E53}.
\end{rem}
L'interprétation géométrique des opérations de transfert peut se faire comme suit : soit $X$ (resp. $Y$) un $K(G,1)$ (resp. un $K(H,1)$) et supposons pour simplifier que $X$ et $Y$ ont une topologie ``raisonnable''. L'inclusion $H\hookrightarrow G$ correspond à un revêtement fini $p:Y\To X$ ; le transfert
$$V_{H\to G}:H^{\bullet}(Y,\CC)\simeq H^{\bullet}(H,\CC)\To H^{\bullet}(G,\CC)\simeq H^{\bullet}(X,\CC)$$
n'est autre que le morphisme d'intégration dans les fibres (ou morphisme de Gysin)
$$p_*: H^{\bullet}(Y,\CC)\To H^{\bullet}(X,\CC).$$
Avec cette interprétation, on retrouve bien le fait mentionné ci-dessus, à savoir
$$V\circ i^*=p_*\circ p^*=\mathrm{deg}(p)\mathrm{Id}=[G:H]\mathrm{Id}.$$

\subsection{Critère de 1-formalité}

Dans ce paragraphe, nous rappelons la notion de 1-formalité d'une algèbre différentielle graduée (\textsc{adg} dans la suite). Pour une discussion plus complète de cette notion, nous renvoyons à \cite{GM}.
\begin{defi}\label{modele minimal}
Une \textsc{adg} $(\mathcal{M},d)$ est dite \emph{1-minimale} si
\begin{enumerate}[(i)]
\item elle est connexe
\item $\mathcal{M}$ peut s'écrire comme une suite d'extension élémentaire (dite de Hirsch)
$$\CC=\mathcal{M}_0\subset\mathcal{M}_1\subset\mathcal{M}_2\subset\dots$$
c'est-à-dire $\mathcal{M}_{i+1}\simeq \mathcal{M}_i\otimes\bigwedge(V_i)$ où $V_i$ est placé en degré 1.
\item $d$ est \emph{décomposable} : au cours de chaque extension élémentaire, $d$ envoie $V_i$ dans $\mathcal{M}_i^+\wedge\mathcal{M}_i^+$, $\mathcal{M}_i^+$ désignant les éléments de degré positif de $\mathcal{M}_i$.
\end{enumerate}

Un morphisme (d'\adg) $\rho:\mathcal{M}\To \mathcal{A}$ est un \emph{1-modèle minimal} pour $(\mathcal{A},d)$ si $(\mathcal{M},d)$ est 1-minimale et si $\rho^*:H^*(\mathcal{M})\To H^*(\mathcal{A})$ induit un isomorphisme en degré 0 et 1 et est injectif en degré 2.
\end{defi}
Un des intérêts de cette définition réside dans la proposition suivante.
\begin{prop}[Sullivan, voir \cite{DGMS}]\label{existence modele min}
Toute \textsc{adg} admet (à isomorphisme près) un unique\footnote{parler d'unicité nécessite d'introduire les notions de points bases et d'homotopies entre \textsc{adg} pour lesquelles nous renvoyons une fois encore à \cite{GM}.} 1-modèle minimal.
\end{prop}
Explicitons la construction dans le cas qui va nous intéresser, à savoir celui de l'algèbre de De Rham $\mathcal{E}^\bullet(X)$ d'une variété différentiable $X$. On souhaite construire inductivement un 1-modèle minimal $M_X^{(1)}\To\mathcal{E}^\bullet(X)$ qui donne un isomorphisme en degré 0 et 1 et un morphisme injectif en degré 2. On commence donc par poser :
$$M_X^{(1)}(1)=\bigwedge(H^1(X,\CC))$$
muni de la différentielle $d_1$ nulle, le morphisme
$$\rho_1:M_X^{(1)}(1)\To \mathcal{E}^\bullet(X)$$
correspondant à un choix de représentants des classes de $H^1(X,\CC)$ fixé une fois pour toute. On a bien un isomorphisme en degré 0 et 1 mais, en degré 2, on a :
$$\rho_1^*:H^2(M_X^{(1)}(1))=\bigwedge^2 H^1(X,\CC)\To H^2(X,\CC).$$
On pose donc
$$V_2=\mathrm{Ker}(\rho_1^*)=\mathrm{Ker}(\bigwedge^2 H^1(X,\CC)\To H^2(X,\CC))$$
et on considère l'extension
$$M_X^{(1)}(2)=M_X^{(1)}(1)\otimes \bigwedge(V_2)$$
et $d_2$ est définie sur $V_2$ comme l'injection naturelle $V_2\hookrightarrow \bigwedge(H^1(X,\CC))=M_X^{(1)}(1)$. Pour définir
$$\rho_2:M_X^{(1)}(2)\To \mathcal{E}^\bullet(X),$$
on le définit sur $V_2$. Un élément $v$ de $V_2$, vu comme 2-classe, est exacte : $v=du$. On pose alors $\rho_2(v)=u$ (à nouveau en faisant un choix de primitive). Examinons l'effet de cette extension au niveau cohomologique. Comme les éléments de $V_2$ ne sont pas fermés (pour $d_2$), on ne change pas la cohomologie en degré 1. En degré 2, on a supprimé le défaut d'injectivité provenant du noyau de
$$\bigwedge^2 H^1(X,\CC)\To H^2(X,\CC)$$
mais on a éventuellement introduit de nouveaux éléments de
$$V_3=\mathrm{Ker}\left(\rho_2^*:H^2(M_X^{(1)}(2))\To H^2(X,\CC)\right).$$
La construction se produit donc inductivement en posant
\begin{align*}
M_X^{(1)}(i+1)&=M_X^{(1)}(i)\otimes \bigwedge(V_{i+1})\quad\textrm{avec}\\
V_{i+1}&=\mathrm{Ker}\left(\rho_i^*:H^2(M_X^{(1)}(i))\To H^2(X,\CC)\right)
\end{align*}
et en construisant $d_{i+1}$ et $\rho_{i+1}$ comme nous l'avons fait pour passer de $M_X^{(1)}(1)$ à $M_X^{(1)}(2)$. Le 1-modèle minimal de $\mathcal{E}^\bullet(X)$ est alors obtenue en prenant la limite inductive de cette suite d'extension :
$$M_X^{(1)}=\bigcup_{i\ge1}M_X^{(1)}(i).$$

\noindent Ceci mène naturellement à la définition suivante.
\begin{defi}\label{1 formalité}
Une variété différentiable $X$ est dite \emph{1-formelle} si son algèbre de De Rahm $(\mathcal{E}^\bullet(X),d)$ l'est, c'est-à-dire si $(\mathcal{E}^\bullet(X),d)$ et $(H^\bullet(X),d)$ ont même 1-modèle minimal.
\end{defi}

Pour finir, signalons le critère suivant de 1-formalité (dû à Morgan) portant uniquement sur le groupe fondamental.
\begin{theo}[th. 9.4, p. 198 \cite{M78}]\label{quadraticité}
Une variété différentiable $X$ (dont le groupe fondamental est de présentation finie) est 1-formelle si et seulement si l'algèbre de Lie $\lie(\gf{X})$ est de présentation quadratique. Ceci est également équivalent à la surjectivité de l'application :
$$H^2(\lie_1(\gf{X}),\RR)\To H^2(\lie(\gf{X}),\RR).$$
\end{theo}
En effet, la correspondance existante entre 1-modèle minimal de $\mathcal{E}^\bullet(X)$ et complétion de Mal\v{c}ev du groupe fondamental est une simple dualité.
\begin{theo}[Sullivan, voir cependant \cite{DGMS}]\label{dualité modele minimal/lie}
Soit $X$ une variété différentiable dont le groupe fondamental est de présentation finie. Le 1-modèle minimal de l'algèbre de De Rham $\mathcal{E}^\bullet(X)$
$$M_X^{(1)}(1)\subset\dots\subset M_X^{(1)}(i)\subset M_X^{(1)}(i+1)\dots$$
et la complétion de Mal\v{c}ev de $\gf{X}$
$$\dots\To \lie_{i+1}(\gf{X})\To \lie_i(\gf{X})\To\dots\To \lie_1(\gf{X})\To 0$$
sont duaux l'un de l'autre.
\end{theo}
\begin{rem}
Le fait que la différentielle vérifie $d^2=0$ se traduit exactement par l'identité de Jacobi au niveau du dual.
\end{rem}

\section{Cas du degré 2}

\subsection{Formalité et groupes nilpotents}

Les résultats de la section précédente s'applique pleinement à la catégorie des variétés kählériennes compactes comme le montre le résultat suivant (dont la démonstration est une conséquence directe du lemme du $dd^c$).
\begin{theo}[\cite{DGMS}]\label{formalité}
Toute variété kählérienne compacte $X$ est formelle ; plus précisément, les algèbres $(\mathcal{E}^\bullet(X),d)$ et $(H^\bullet(X),0)$ sont équivalentes \emph{via} l'algèbre $(\mathcal{E}^{\bullet}_{d^c}(X),d)$ des formes $d^c$-fermées. En particulier, une variété kählérienne compacte est 1-formelle.
\end{theo}
Rappelons la
\begin{defi}
Une \textsc{adg} $(\mathcal{A},d)$ est dite \emph{formelle} si elle est équivalente à sa propre algèbre de cohomologie (avec différentielle nulle) ; c'est-à-dire si il existe une chaîne de quasi-isomorphismes :
$$(\mathcal{A},d)\longleftarrow (\mathcal{C}_1,d_1)\To (\mathcal{C}_2,d_2)\longleftarrow\dots\longleftarrow(\mathcal{C}_n,d_n)\To (H^\bullet(\mathcal{A}),0).$$
Une variété différentiable $X$ est dite \emph{formelle} si son algèbre de De Rham $(\mathcal{E}^\bullet(X),d)$ l'est.
\end{defi}

Nous pouvons dès à présent démontrer le théorème \ref{th2} grâce au critère de quadraticité.\\

\begin{demo}[du théorème \ref{th2}]
Soit $X$ une variété kählérienne compacte dont le groupe fondamental est nilpotent sans torsion. Le théorème \ref{cohomologie nilp} s'applique et on a :
$$\forall k\ge0,\,H^k(\gf{X},\RR)\simeq H^k(\lie(\gf{X}),\RR).$$
Or, d'après les théorèmes \ref{quadraticité} et \ref{formalité}, on sait que la flèche naturelle
$$H^2(\gf{X}_{ab},\RR)\simeq H^2(\lie_1(\gf{X}),\RR)\To H^2(\lie(\gf{X}),\RR)\simeq H^2(\gf{X},\RR)$$
est surjective. Or, comme le groupe de gauche est aussi
$$H^2(\gf{X}_{ab},\RR)\simeq H^1(X,\RR)\bigwedge H^1(X,\RR),$$
on a bien :
$$H^2(\gf{X},\RR)=\mathrm{Im}\left(H^1(X,\RR)\bigwedge H^1(X,\RR)\To H^2(X,\RR)\right).$$

Dans le cas général ($\gf{X}$ virtuellement nilpotent), on sait que $X$ admet un revêtement galoisien fini $Y\To X$ de groupe de Galois $G=\gf{X}/\gf{Y}$ et tel que $\gf{Y}$ est nilpotent sans torsion. On peut appliquer la discussion précédente à $Y$ et on obtient donc un morphisme surjectif :
$$H^2(\gf{Y}_{ab},\RR)\To H^2(\gf{Y},\RR).$$
En considérant les éléments $G$-invariants de ces deux espaces ($G$ agit naturellement sur $\gf{Y}_{ab}$) et en appliquant la proposition \ref{transfert}, on obtient la même conclusion pour $H^2(\gf{X},\RR)$.

Tout ceci montre en particulier que $H^2(\gf{X},\CC)$ est une sous-structure de Hodge de $H^2(X,\CC)$ (comme image d'un morphisme de structure de Hodge). Enfin, le fait que
$$\mathrm{Im}\left(H^1(X,\CC)\bigwedge H^1(X,\CC)\To H^2(X,\CC)\right)\neq0$$
(pour une variété kählérienne compacte) est une conséquence directe du théorème de Lefschetz difficile et du fait que $H^1(X,\CC)$ est lui-même non nul (un groupe nilpotent infini admet des quotients abéliens infinis).
\end{demo}

\begin{rem}
La démonstration ci-dessus ne nécessite en réalité que la 1-formalité (et le caractère nilpotent du groupe fondamental) de $X$.
\end{rem}

\subsection{Exemples non-kählériens}

Nous venons de montrer dans la section précédente que, pour un groupe kählérien nilpotent $\GG$, l'application naturelle
$$H^2(\GG_{ab},\CC)\To H^2(\GG,\CC)$$
était surjective. Pour nous convaincre que ceci est bien spécifique au cas kählérien (au moins au cas des variétés 1-formelles), voici quelques exemples.

\begin{ex}\label{exemple1}
Soit $G$ le groupe de Heisenberg réel, $\GG$ le réseau des éléments de $G$ à coefficients dans $\ZZ$ et considérons la variété différentiable $X=G/\GG$. Le groupe $\GG$ s'écrit donc comme une extension centrale
$$1\To \ZZ\To \GG\To \GG_{ab}\simeq \ZZ^2\To 1.$$
Cette décomposition induit une suite exacte
$$0\To H^1(\ZZ)\To H^2(\GG_{ab})\To H^2(\GG).$$
Pour des raisons de dimension, on constate immédiatement que la flèche
$$H^2(\GG_{ab})\To H^2(\GG)$$
est identiquement nulle alors que $H^2(\GG)$ ne l'est pas. En effet, comme $X$ est un $K(\GG,1)$, la dualité de Poincaré entraîne :
$$H^2(\GG)\simeq H^2(X) \simeq H^1(X)^*\simeq H^1(\GG)^*\neq 0.$$
\end{ex}

\begin{ex}\label{exemple2}
Pour obtenir un exemple de variété complexe (plus proche de la situation kählérienne), on reprend l'exemple précedent mais avec cette fois des coefficients complexes. Soit donc $G$ le groupe de Heisenberg complexe, $\GG$ le réseau des éléments de $G$ à coefficients dans $\ZZ[i]$ et considérons la variété complexe (non-kählérienne) $X=G/\GG$. Comme $\GG_{ab}\simeq\ZZ^4$, le groupe $H^2(\GG_{ab})$ est de dimension 6. Or, $X$ est à nouveau un $K(\GG,1)$ et on a donc $H^2(\GG)\simeq H^2(X)$ et il est bien connu que $b_2(X)=8$. La flèche $H^2(\GG_{ab})\To H^2(\GG)$ ne peut donc pas être surjective.
\end{ex}

\section{Structure de Hodge en degrés supérieurs}

\subsection{La \textsc{shm} de Morgan}

Le théorème \ref{th1} est en fait une réécriture des résultats de Morgan \cite{M78}. En effet, d'après \cite{M78}, on peut munir le 1-modèle minimal d'une variété kählérienne compacte d'une structure de Hodge mixte (\textsc{shm}\footnote{pour les notions concernant les \textsc{shm}, nous renvoyons à \cite{PS}.} dans la suite) fonctorielle. Plus précisément, si $X$ est une variété kählérienne compacte, notons $\rho_X:M^{(1)}_X\To\mathcal{E}^\bullet(X)$ le 1-modèle minimal de son algèbre de De Rham. Comme $X$ est formelle (théorème \ref{formalité}), $(\mathcal{E}^\bullet(X),d)$ et $(H^\bullet(X),0)$ ont même 1-modèle minimal et on dispose d'un morphisme
$$\sigma_X:M^{(1)}_X\To H^\bullet(X).$$
\begin{theo}\label{morgan}
Avec les notations ci-dessus, l'algèbre $(M^{(1)}_X,d^{(1)})$ possède une \textsc{shm} fonctorielle vérifiant :
\begin{enumerate}[(1)]
\item la différentielle $d^{(1)}$ et le produit dans l'algèbre $M^{(1)}_X$ sont des morphismes de \textsc{shm}.
\item l'application $\sigma_X:M^{(1)}_X\To H^\bullet(X,\CC)$ est un morphisme de \textsc{shm}.
\end{enumerate}
\end{theo}
La filtration par le poids $W_\bullet$ de $M^{(1)}_X$ est donnée par la description de $M^{(1)}_X$ comme l'union croissante des sous-algèbres $M^{(1)}_X(n)$ (et est donc duale de la suite centrale descendante, la dualité étant fournie par le théorème \ref{dualité modele minimal/lie}). La filtration de Hodge $F^\bullet$ provient elle de celle de $H^1(X,\CC)$.\\

Comme mentionné ci-dessus, le théorème \ref{th1} consiste maintenant à réinterpréter les résultats de Morgan en termes de cohomologie du groupe $\gf{X}$ (dans le cas nilpotent).\\
\begin{demo}[du théorème \ref{th1}]
Soit donc $X$ dont le groupe fondamental est (dans un premier temps) nilpotent sans torsion. D'après les théorèmes \ref{coho lie} et \ref{dualité modele minimal/lie}, on dispose des isomorphismes :
$$H^*(\gf{X},\CC)\simeq H^*(\lie(\gf{X}),\CC)\simeq H^*(M^{(1)}_X).$$
On peut alors appliquer le théorème \ref{morgan} ; comme $M^{(1)}_X$ admet une \textsc{shm} qui fait de la différentielle un morphisme de \textsc{shm}, cette structure passe en cohomologie et ce procédé nous permet donc de définir une \textsc{shm} sur la cohomologie de $\gf{X}$. D'autre part, comme le morphisme
$$\sigma_X:M^{(1)}_X\To H^*(X,\CC)$$
est à la fois un morphisme d'algèbres et un morphisme de \textsc{shm}, le morphisme induit
$$\sigma_X^*:H^*(\gf{X},\CC)\simeq H^*(M^{(1)}_X)\To H^*(X,\CC)$$
est bien un morphisme de \textsc{shm}.

Si $\gf{X}$ est seulement supposé virtuellement nilpotent, on sait qu'il admet un sous-groupe d'indice fini nilpotent sans torsion. Si on note $Y\To X$ le revêtement étale fini (galoisien) correspondant à ce sous-groupe, la discussion ci-dessus s'applique à $\gf{Y}$ et on peut donc munir $H^*(\gf{Y},\CC)$ d'une \textsc{shm} fonctorielle. Comme le groupe de Galois $G=\gf{X}/\gf{Y}$ agit par biholomorphismes sur $Y$, l'action de $G$ sur $H^*(\gf{Y},\CC)$ préserve donc la \textsc{shm} et ceci montre que
$$H^*(\gf{X},\CC)=H^*(\gf{Y},\CC)^{G}$$
hérite d'une \textsc{shm} et que le morphisme naturel
$$H^*(\gf{X},\CC)=H^*(\gf{Y},\CC)^{G}\To H^*(X,\CC)=H^*(Y,\CC)^G$$
est bien un morphisme de \textsc{shm}.
\end{demo}

En guise de conclusion, récapitulons les différents isomorphismes (et morphismes) qui nous ont permis de munir $H^*(\GG,\CC)$ d'une \textsc{shm} dans le cas nilpotent sans torsion. Dans le diagramme suivant
$$\xymatrix{H^*(\GG,\CC)\ar[d]\ar[r]^{\!\!\!\!\!\!\!\!\!\!\!\!\!\!\sim}_{\!\!\!\!\!\!\!\!\!\!\!\!\!\!(1)} & H^*(\lie(\gf{X}),\CC)\ar[r]^{\,\,\,\sim}_{\,\,\,(2)} & H^*\left(M^{(1)}_{\mathcal{E}^\bullet(X)}\right) \ar[dl]_{\sim}^{(3)}\\
H^*(X,\CC) & H^*\left(M^{(1)}_X\right) \ar[l]^{(4)}&
}$$
on a (volontairement) noté $M^{(1)}_{\mathcal{E}^\bullet(X)}$ (resp. $M^{(1)}_X$) le 1-modèle minimal de l'algèbre de De Rham de $X$ (resp. celui de l'agèbre de cohomologie $H^*(X)$). L'isomorphisme (1) est donné par le théorème de Nomizu \ref{cohomologie nilp} et la flèche (2) correspond au théorème de Sullivan \ref{dualité modele minimal/lie}. La formalité (théorème \ref{formalité}) quant à elle assure que les modèles minimaux $M^{(1)}_{\mathcal{E}^\bullet(X)}$ et $M^{(1)}_X$ sont isomorphes et fournit la flèche (3). Le théorème de Morgan \ref{morgan} montre enfin comment munir la cohomologie de $M^{(1)}_X$ d'une \textsc{shm} à partir de la structure de Hodge de $H^*(X,\CC)$ de telle sorte que la flèche (4) soit un morphisme de \textsc{shm} et complète le parcours de ce diagramme.


\subsection{Revue des exemples connus}\label{exemple-fred}

Les seuls exemples connus de groupes kählériens nilpotents (non-abéliens) sont ceux exhibés dans \cite{Ca95nilp} et \cite{SV86} et bien sûr leurs produits. Ils sont tous obtenus comme extension centrale de $\ZZ$ par une groupe abélien $A$ (sans torsion, de rang $\ge8$) :
$$1\To \ZZ\To \GG\To A\To 1.$$

Pour la commodité du lecteur, redonnons une des constructions de \cite{Ca95nilp}. Soit $V$ le complémentaire dans $\PP^{2n+1}$ de deux sous-espaces linéaires de dimension $n$ en position générale (\emph{i.e.} deux copies de $\PP^n$ disjointes). Il est bien connu que $V$ admet une structure de $\CC^*$-fibré sur $\PP^n\times\PP^n$. Si $A$ est une variété abélienne admettant une application homolorphe, surjective et finie sur $\PP^n$, considérons $Y$ le $\CC^*$-fibré sur $A\times A$ obtenu par tiré en arrière :
$$\xymatrix{Y\ar[r]^{g}\ar[d] & V\ar[d]\\
A\times A\ar[r]^f & \PP^n\times\PP^n
}.$$
La structure de $\CC^*$-fibré sur $A\times A$ montre que le groupe fondamental de $Y$ est une extension
$$1\To \ZZ=\gf{\CC^*}\To \gf{Y}\To \gf{A}\times\gf{A}\To 1$$
qui est en fait centrale et que $Y$ est un $K(\gf{Y},1)$. La variété $Y$ est quasi-projective mais on peut également réaliser $\gf{Y}$ comme le groupe fondamental d'une variété \emph{projective} ; il faut pour cela utiliser les résultats de théorie de Morse stratifiée de Goreski et MacPherson. En effet, si $L$ désigne un sous-espace linéaire de dimension $n$ contenu dans $V$ et en position générale, on peut appliquer les résultats de \cite[th. p. 195]{GMc} :
$$\forall\, i\leq n-1,\, \pi_i(Y,h^{-1}(L))=0.$$
Si $n\ge 2$, on obtient en particulier en posant $X=h^{-1}(L)$ : $\gf{X}=\gf{Y}$. Pour $n\ge2$, le groupe $\gf{Y}$ est donc aussi le groupe fondamental de la variété projective $X$ (qui est lisse pour un choix de $L$ générique). Nous allons voir que l'on peut vérifier \emph{à la main} que la cohomologie du groupe $\gf{X}$ satisfait aux conclusions des théorèmes \ref{th1} et \ref{th2}.\\

\noindent\textbf{Première méthode :}

Comme $Y$ est un $K(\gf{X},1)$, on sait que $H^*(\gf{X},\CC)\simeq H^*(Y,\CC)$ et l'injection canonique $j:X\hookrightarrow Y$ induit les morphismes :
$$H^*(\gf{X},\CC)\simeq H^*(Y,\CC)\stackrel{j^*}{\To}H^*(X,\CC).$$
Or, d'après \cite{D71}, la cohomologie de $Y$ (qui est quasi-projective lisse) porte une structure de Hodge mixte et on sait également que le morphisme $j$ induit un morphisme de structure de Hodge mixte en cohomologie. C'est exactement ce qui est prédit par le théorème \ref{th1}.$\square$\\

\noindent\textbf{Deuxième méthode :}

Le groupe fondamental $\GG=\gf{X}$ s'écrit comme une extension centrale :
$$1\To \ZZ\To \GG\To \GG_{ab}=\gf{A}\times\gf{A}\To 1\quad (*),$$
et celle-ci permet de calculer la cohomologie de $\GG$ grâce à la suite spectrale de Hochschild-Serre \cite{HS}. Cette suite spectrale ayant peu de termes non nuls, elle dégénère en $E_3$ et induit une suite longue :
$$\xymatrix{E^{2,0}_2\ar@{=}[d] \ar[r] & H^2(\GG) \ar@{=}[d] \ar[r] & E^{1,1}_2\ar@{=}[d] \ar[r]^{d_2} & E^{3,0}_2\ar@{=}[d] \ar[r] & \dots\\
H^2(\GG_{ab})\ar[r]& H^2(\GG)\ar[r]& H^1(\GG_{ab})\ar[r]^{cl}& H^3(\GG_{ab})\ar[r]& \dots}$$
Le morphisme $cl$ ci-dessus est donné par le produit par la classe de l'extension $(*)$, qui est aussi celle du $\CC^*$-fibré $Y\To A\times A$. Si $f:A\To\PP^n$ désigne la projection, la classe d'extension de $(*)$ est donnée par
$$cl=(f^*\omega,-f^*\omega)\in H^2(\GG_{ab})=H^2(A\times A,\CC)$$
où $\omega$ désigne la classe hyperplane de $\PP^n$. Cette classe est donc de type (1,1) et non-dégénérée. On constate alors facilement que le produit par cette classe est injectif sur le $H^1(\GG_{ab})$ ; de façon équivalente, la flèche $H^2(\GG_{ab})\To H^2(\GG)$ est surjective (c'est le contenu du théorème \ref{th2}). De plus, la suite longue ci-dessus montre que les groupes de cohomologie de $\GG$ se décomposent de la façon suivante :
\begin{align}\label{shm explicite}
0\To H^k(\GG_{ab})/\left(H^{k-2}(\GG_{ab})\wedge cl\right) &\To H^k(\GG,\CC) \To\\ &\mathrm{Ker}\left(H^{k-1}(\GG_{ab})\stackrel{cl}{\To}H^{k+1}(\GG_{ab})\right)\To 0\nonumber.
\end{align}
Comme la classe $cl$ est de type (1,1) et que les groupes $H^j(\GG_{ab})$ sont naturellement munis de structure de Hodge (pures), la décomposition (\ref{shm explicite}) donne une description de la \textsc{shm} sur $H^k(\GG,\CC)$ (au moins des structures de Hodge des quotients successifs de la filtration par le poids).
$\square$\\

\begin{rem}
Les morphismes $H^k(\GG)\To H^{k-1}(\GG_{ab})$ obtenus à partir de la suite spectrale ci-dessus peuvent être exprimés explicitement. Par exemple, pour $k=2$, si $f\in C^2(\GG,\CC)$ est un cocyle et si $z\in\ZZ$ est un générateur du centre $\ZZ$ de $\GG$, l'application
$$[f]:x\mapsto f(x,z)-f(z,x)$$
définit un élément de $H^1(\GG_{ab})$.
\end{rem}

Pour finir, notons que la discussion ci-dessus s'appliquent pour les exemples de \cite{SV86} ; en effet, les groupes kählériens obtenus sont encore des extensions centrales
$$1\To \ZZ\To \gf{X}\To Q \To 1$$
où $Q$ est un groupe abélien (le groupe fondamental d'une variété abélienne).

\bibliographystyle{smfalpha}
\bibliography{myref}

\providecommand{\bysame}{\leavevmode ---\ }
\providecommand{\og}{``}
\providecommand{\fg}{''}
\providecommand{\smfandname}{\&}
\providecommand{\smfedsname}{\'eds.}
\providecommand{\smfedname}{\'ed.}
\providecommand{\smfmastersthesisname}{M\'emoire}
\providecommand{\smfphdthesisname}{Th\`ese}
\begin{thebibliography}{DGMS75}

\bibitem[Bro82]{B87}
{\scshape K.~S. Brown} -- \emph{Cohomology of groups}, Graduate Texts in
  Mathematics, vol.~87, Springer-Verlag, New York, 1982.

\bibitem[Cam95]{Ca95nilp}
{\scshape F.~Campana} -- {\og Remarques sur les groupes de {K}\"ahler
  nilpotents\fg}, \emph{Ann. Sci. \'Ecole Norm. Sup. (4)} \textbf{28} (1995),
  no.~3, p.~307--316.

\bibitem[Del71]{D71}
{\scshape P.~Deligne} -- {\og Th\'eorie de {H}odge. {II}\fg}, \emph{Inst.
  Hautes \'Etudes Sci. Publ. Math.} (1971), no.~40, p.~5--57.

\bibitem[Del06]{D}
{\scshape T.~Delzant} -- {\og L'invariant de {B}ieri {N}eumann {S}trebel des
  groupes fondamentaux des vari{\'e}t{\'e}s k{\"a}hl{\'e}riennes\fg}, preprint
  arXiv:math/0603038, {\`a} para{\^i}tre dans Math. Annalen, 2006.

\bibitem[DGMS75]{DGMS}
{\scshape P.~Deligne, P.~Griffiths, J.~Morgan {\normalfont \smfandname}
  D.~Sullivan} -- {\og Real homotopy theory of {K}\"ahler manifolds\fg},
  \emph{Invent. Math.} \textbf{29} (1975), no.~3, p.~245--274.

\bibitem[DPS09]{Dim09}
{\scshape A.~Dimca, c.~Papadima {\normalfont \smfandname} A.~Suciu} -- {\og
  Non-finiteness properties of the fundamental groups of smooth projective
  varieties\fg}, \emph{J. Reine und Angew. Math.} \textbf{629} (2009),
  p.~89--105.

\bibitem[Eck53]{E53}
{\scshape B.~Eckmann} -- {\og Cohomology of groups and transfer\fg}, \emph{Ann.
  of Math. (2)} \textbf{58} (1953), p.~481--493.

\bibitem[GM81]{GM}
{\scshape P.~A. Griffiths {\normalfont \smfandname} J.~W. Morgan} --
  \emph{Rational homotopy theory and differential forms}, Progress in
  Mathematics, vol.~16, Birkh\"auser Boston, Mass., 1981.

\bibitem[GM88]{GMc}
{\scshape M.~Goresky {\normalfont \smfandname} R.~MacPherson} --
  \emph{Stratified {M}orse theory}, Ergebnisse der Mathematik und ihrer
  Grenzgebiete (3) [Results in Mathematics and Related Areas (3)], vol.~14,
  Springer-Verlag, Berlin, 1988.

\bibitem[HS53]{HS}
{\scshape G.~Hochschild {\normalfont \smfandname} J.-P. Serre} -- {\og
  Cohomology of group extensions\fg}, \emph{Trans. Amer. Math. Soc.}
  \textbf{74} (1953), p.~110--134.

\bibitem[Kol95]{K95}
{\scshape J.~Koll{\'a}r} -- \emph{Shafarevich maps and automorphic forms}, M.
  B. Porter Lectures, Princeton University Press, Princeton, NJ, 1995.

\bibitem[Mal49]{M49}
{\scshape A.~I. Mal\v{c}ev} -- {\og On a class of homogeneous spaces\fg},
  \emph{Izvestiya Akad. Nauk. SSSR. Ser. Mat.} \textbf{13} (1949), p.~9--32.

\bibitem[Mor78]{M78}
{\scshape J.~W. Morgan} -- {\og The algebraic topology of smooth algebraic
  varieties\fg}, \emph{Inst. Hautes \'Etudes Sci. Publ. Math.} (1978), no.~48,
  p.~137--204.

\bibitem[Nom54]{N54}
{\scshape K.~Nomizu} -- {\og On the cohomology of compact homogeneous spaces of
  nilpotent {L}ie groups\fg}, \emph{Ann. of Math. (2)} \textbf{59} (1954),
  p.~531--538.

\bibitem[PS08]{PS}
{\scshape C.~A.~M. Peters {\normalfont \smfandname} J.~H.~M. Steenbrink} --
  \emph{Mixed {H}odge structures}, Ergebnisse der Mathematik und ihrer
  Grenzgebiete. 3. Folge. A Series of Modern Surveys in Mathematics, vol.~52,
  Springer-Verlag, Berlin, 2008.

\bibitem[Rag72]{R72}
{\scshape M.~S. Raghunathan} -- \emph{Discrete subgroups of {L}ie groups},
  Springer-Verlag, New York, 1972, Ergebnisse der Mathematik und ihrer
  Grenzgebiete, Band 68.

\bibitem[Sta63]{St63}
{\scshape J.~Stallings} -- {\og A finitely presented group whose 3-dimensional
  integral homology is not finitely generated\fg}, \emph{Amer. J. Math.}
  \textbf{85} (1963), p.~541--543.

\bibitem[SVdV86]{SV86}
{\scshape A.~J. Sommese {\normalfont \smfandname} A.~Van~de Ven} -- {\og
  Homotopy groups of pullbacks of varieties\fg}, \emph{Nagoya Math. J.}
  \textbf{102} (1986), p.~79--90.

\end{thebibliography}

\vspace*{0.5cm}
\begin{flushright}
\begin{minipage}{5cm}
Beno\^it \textsc{Claudon}\\
Institut Fourier - UMR 5582\\
100, rue des Maths\\
B.P. 74\\
38402 Saint-Martin d'Hères\\
France

\vspace*{0.3cm}
\noindent Benoit.Claudon@ujf-grenoble.fr
\end{minipage}
\end{flushright}

\end{document}